\documentclass{amsart}      
\usepackage{amssymb}      
\usepackage{latexsym}      
     
\newcommand{\R}{{\mathbb R}}

\newtheorem{lemma}{Lemma}[section]

\begin{document}         
\title[Aperiodic linearly repetitive Delone sets]{Aperiodic linearly
  repetitive Delone sets are densely repetitive}         
\author[D.~Lenz]{Daniel  
  Lenz$\,^{1}$}    
   
\maketitle         
\vspace{0.3cm}             
\noindent             
$^1$ Fakult\"at f\"ur Mathematik, TU Chemnitz, D-09107 Chemnitz, Germany, E-Mail: dlenz@mathematik.tu-chemnitz.de\\[0.3cm]             
         
\begin{abstract}   We show that aperiodic linearly repetitive Delone
  sets are densely repetitive. This confirms a conjecture of Lagarias and Pleasants.
\end{abstract}         
         
\section{Introduction}  \label{Introduction}      
In \cite{LP}, Lagarias and Pleasants study the problem of characterizing
the simplest aperiodic discrete point sets. To this aim they carry out a careful study of
linearly repetitive and densely repetitive Delone sets. As for the
relationship between these two concepts they formulate the following
conjecture:

\medskip

{\bf Conjecture} (= Conjecture 1.2 a in \cite{LP}).  Every aperiodic
linearly repetitive Delone set is densely repetitive. 

\medskip

It is the aim of this note to prove the conjecture. The proof is
a rather direct  consequence of a reformulation of the conjecture in
terms of  lower bounds on local complexity as given already in \cite{LP} and   two
essentially known  facts. The first fact is a
certain repulsion property of occurences of the same patches in
aperiodic linearly repetitive Delone sets (see the work of Solomyak
\cite{Sol} and Durand \cite{Du2}). The
second is a connection between absence of  ``local  periods`` and
lower bounds on local complexity as studied e.g. by the author in \cite{Len}.

\medskip

The next section recalls the necessary notation (mostly taken from \cite{LP})
and gives a proof of the conjecture. 

\section{Notations and proof}
A set $X\subset \R^d$ is called Delone set if there exist $0< r(X), R(X)
<\infty$ such that every ball in $\R^d$ with radius $r(X)$  
meets at most one point of $X$ and every ball with radius $R(X)$ meets
at least one point of $X$. For $T>0$, a set of the form $X\cap B(x,T)$
with $x\in X$ is called a $T$-patch centered at $x$. Here, $B(p,s)$
denotes the closed ball arround $p$ with radius $s$. The number $N_X (T)$  of different
$T$-patches up to translation is  defined by 
$$ N_X (T) \equiv \sharp \{ (X- x)\cap B(0,T) : x\in X\},$$
where $\sharp S$ denotes the number of elements of the set $S$. The
Delone set $X$ is called repetitive if, for every $T>0$, there
exists a  finite number $M$ such that every closed ball of radius $M$ in
$\R^d$ contains the center of a translate of every possible $T$-patch in
$X$. The smallest such $M $ is denoted by $M_X (T)$. 

\medskip

A Delone set $X$ is
called linearly repetitive if there exists a constant  $C_{LR} (X)\in \R$ with
$M_X (T) \leq C_{LR} (X)\, T$ for every $T\geq 1$. It is called densely
repetitive if there exists a constant $C_{DR} (X)\in \R$ with   $M_X (T)
\leq  C_{DR} (X)  N_X^{\frac{1}{d}}(T)$ for $T \geq 1$.  The set $X$  is called
non-periodic  if  it is not equal to a translate of itself 
and it is called aperiodic if this property holds
for all elements in the hull. Here, the hull is the closure of its
translates in the natural topology \cite{LP}. For repetitive Delone sets
non-periodicity and aperiodicity are equivalent. 

\medskip

Our proof of the conjecture is based on the following two lemmas. The
first gives a precise version of the repulsion property mentioned in the
introduction. 

\begin{lemma} \label{repulsion} Let $X$ be an aperiodic linearly
  repetitive Delone set. Then, there exists a constant $\kappa(X)>0$
  such that $\| x -y \| \geq \kappa (X) T$ whenever $x,y\in X$ with
  $x\neq y $ and  $(X-x)\cap B(x,T) = (X-y)  \cap B(y,T)$ for some
  $T>1$. Here, $\|\cdot\|$ denotes the Euclidean norm. 
\end{lemma}

For linearly repetitive tilings this lemma  (and in fact a slightly
stronger version) is proven in Lemma 2.4 of  \cite{Sol}  (see  proof of
Theorem 2.2 in  \cite{LP2}  and \cite{Du2}  as well). The
proof given there is formulated in terms of tilings associated to
primitive substitutions but only uses linear repetitivity (called strong
repetitivity in \cite{Sol}). It can easily be carried over to Delone
sets.  For completeness reasons, we include a short discussion giving a
{\it proof of Lemma \ref{repulsion}}: Assume the contrary. Then, there
exist $x,y\in X$, $x\neq y$ with 
$$ (X- x) \cap B(x,T) = (X-y) \cap B(y,T)$$
and 
$$ r(X) \leq \| x- y\| \leq  \frac{T}{( C_{LR} +1) ( r(X)^{-1} + 1) }.$$
In particular, for every point $s\in X \cap B(x,T)$, the point $s +
(y-x)$ belongs to $X$ as well. By linear repetitivity, $X \cap B(x,T)$
contains the center of a translate of every $ C_{LR}^{-1} T$-patch  of
$X$. As $C_{LR}^{-1} T > \|x-y\|$, we infer that for every $z\in X$ the
point $z + (y-x)$ belongs to $X$ as well. This proves Lemma \ref{repulsion}. 

\medskip

The second lemma is a consequence of the preceeding lemma. 

\begin{lemma}\label{complexity} Let $X$ be an aperiodic linearly
  repetitive Delone set. Then, there exist constants $\lambda >0$  and
  $T_0 >0$ with $N_X (T)  \geq \lambda T^d$ for $T\geq T_0$. 
\end{lemma}
{\it Proof.} As $X$ is a Delone set, there exists $\lambda_1 >0$ and $T_1
>0$ with 
\begin{equation} \label{anzahl}
\sharp X \cap B(p,T) \geq \lambda_1 T^d
\end{equation} 
for all $p\in R^d$ and $T\geq T_1$. Now, let $\kappa (X)$ be as given in
the previous lemma and consider for $T>1$ the $T$-patches
$$ (X-x) \cap B(0,T)  \;\: \mbox{with}\;\: x\in X \cap B( 0,
\frac{\kappa(X) T }{3}). $$
Then, by the previous lemma, these $T$-patches are pairwise
different. Thus,  using \eqref{anzahl}, we can calculate

\begin{eqnarray*}
N_X (T)
& \geq &\sharp (X \cap  B( 0, \frac{\kappa(X) T }{3}))
\geq \lambda_1  \left( \frac{\kappa(X)T}{3}
\right)^d =   \lambda T^d, 
\end{eqnarray*} 
for $T\geq T_0 \equiv  3 \kappa(X)^{-1} T_1$ and $\lambda = \lambda_1
(\kappa(X) 3^{-1})^d$. \hfill \qedsymbol

\medskip

We can now provide the {\it Proof of the conjecture}: 
It suffices to show
\begin{equation}\label{equivalent}
\liminf_{T\to \infty} \frac{ N_X (T)}{T^d} >0.
\end{equation}
In fact, by the discussion
in  Section 8 of Lagarias/Pleasants \cite{LP} (see \cite{LP2} as well),
validity of \eqref{equivalent} is even equivalent to validity of the
conjecture. Now, \eqref{equivalent}  follows  immediately from
the previous lemma. The conjecture is proven. \hfill \qedsymbol

\medskip

{\bf Acknowledgments.} The author would like to thank Peter Pleasants
and Boris Solomyak for useful discussions.

\end{document}